\documentclass[12pt]{article}
\usepackage{amsfonts}
\usepackage{amssymb}
\usepackage{amsthm}

\newtheorem{theorem}{Theorem}[section]
\newtheorem{lemma}{Lemma}[section]
\newtheorem*{remark}{Remark}
\newtheorem*{corollary}{Corollary}

\author{Yuri A. Turygin\thanks{Partially supported by the Russian Foundation
for Basic Research Grant No. 02-01-00014}}
\title{Approximation of $k$-dimensional maps}
\date{}

\begin{document}
\maketitle

\begin{abstract}
\scriptsize
In this paper we prove the equivalence of the questions of
B.A. Pasynkov (\cite{pas}) and V.V. Uspenskij (\cite{usp}).
We also get some partial results answering these
questions in affirmative. As a corollary to these results we get
an extention of the Hurewicz formula to the extensional dimension.
\end{abstract}

\section {Introduction}

All spaces are assumed to be separable metrizable. By a map we
mean a continuous function, ${\bf I}=[0;1].$ If ${\cal K}$ is a
simplicial complex then by $|{\cal K}|$ we mean the corresponding
polyhedron. By a simplicial map we mean a map $f\colon |{\cal
K}|\longrightarrow |{\cal L}|$ which sends simplices to simplices
and is affine on them. We say that a map $f\colon X\longrightarrow
Y$ has dimension at most $k$ ($\dim f\le k$) if and only if the
dimension of each of its fibers is at most $k$. We recall that a
space $X$ is a $C$-space or has property $C$ if for any sequence
$\{\alpha_n :n\in \mathbb N\}$ of open covers of $X$ there exists
a sequence $\{\mu_n :n\in \mathbb N\}$ of disjoint families of
open sets such that each $\mu_n$ refines $\alpha_n $ and the
union of all systems $\mu_n$ is a cover of $X$. Each
finite-dimensional paracompact space and each
countable-dimensional metrizable space has property $C$. By a
$C$-compactum we mean a compact $C$-space.

In \cite{usp} V.V. Uspenskij introduced the notion of a map
admitting an approximation by $k$-dimensional simplicial maps.
Following him we say that a map $f\colon X \longrightarrow Y$
admits approximation by $k$-dimensional simplicial maps if for
every pair of open covers $\omega_X$ of the space $X$ and
$\omega_Y$ of the space $Y$ there exists a commutative diagram of
the following form

$$
\begin{array}{rcc}
 X &\stackrel{\kappa_X}{\longrightarrow}& |{\cal K}| \\
 f \downarrow &&\downarrow\lefteqn{p}\\
 Y &\stackrel{\kappa_Y}{\longrightarrow}& |{\cal L}|,
\end{array}
$$
\\
where $\kappa_X$ is an $\omega_X$-map, $\kappa_Y$ is an $\omega_Y$-map
and $p$ is a $k$-dimensional simplicial map between polyhedra
$|\cal K|$ and $|\cal L|$.

In that paper V.V. Uspenskij proposed the following question and
conjectured that in the general case the answer to it is "no".

{\bigskip} ${\bf (Q_1)}$ Does every $k$-dimensional map $f\colon
X\longrightarrow Y$ between compacta admit approximation by
$k$-dimensional simplicial maps? {\bigskip}

In \cite{drusp} A.N. Dranishnikov and V.V. Uspenskij proved that
light maps admit approximation by finite-to-one simplicial maps.
In this paper we give some partial results answering the question
of V.V. Uspenskij in affirmative.

\begin{theorem}
\label{th_C-comp} Let $f\colon X \longrightarrow Y$ be a
$k$-dimensional map between $C$-compacta. Then for any pair of
open covers $\omega_X$ of the space $X$ and $\omega_Y$ of the
space $Y$ there exists a commutative diagram of the following form

$$
\begin{array}{rcc}
 X &\stackrel{\kappa_X}{\longrightarrow}& |{\cal K}| \\
 f \downarrow &&\downarrow\lefteqn{p}\\
 Y &\stackrel{\kappa_Y}{\longrightarrow}& |{\cal L}|,
\end{array}
$$
\\
where $\kappa_X$ is an $\omega_X$-map, $\kappa_Y$
is an $\omega_Y$-map and $p$ is a
$k$-dimensional simplicial map between compact polyhedra
$|\cal K|$ and $|\cal L|.$ Furthermore, one can always assume that
$$\dim |{\cal L}| \leq \dim Y \  \  and\  \  \dim |{\cal K}| \leq \dim Y + k.$$
\end{theorem}

\begin{theorem}
\label{arbcomp}
$k$-dimensional maps between compacta admit
approximation by $(k+1)$-dimensional simplicial maps.
\end{theorem}

\begin{theorem}
\label{bingcomp}
$k$-dimensional maps of Bing compacta (i.e. compacta with each
component hereditarily indecomposable) admit
approximation by \\$k$-dimensional simplicial maps.
\end{theorem}

\begin{remark}
Theorem \ref{th_C-comp} was announced by V.V. Uspenskij at the
$8^{th}$ International Topology Symposium in Prague (1996), but the
proof appeared only here.
\end{remark}

It turned out that the question ${\bf (Q_1)}$ is closely related
to the next question proposed by B.A. Pasynkov in \cite{pas}. We
recall that the {\sl diagonal product} of two maps $f\colon
X\longrightarrow Y$ and $g\colon X\longrightarrow Z$ is a map
$f\bigtriangleup g\colon X\longrightarrow Y\times Z$ defined by
$f\bigtriangleup g (x) = (f(x),g(x))$.

{\bigskip} ${\bf (Q_2)}$ Let $f\colon X \longrightarrow Y$ be a
$k$-dimensional map between compacta. Does there exist a map
$g\colon X \longrightarrow {\bf I}^k$ such that $\dim(f
\bigtriangleup g) \leq 0?$ {\bigskip}

In this paper we prove the following theorem which states that
the questions ${\bf (Q_1)}$ and ${\bf (Q_2)}$ are equivalent.

\begin{theorem}
\label{PU} Let $f\colon X\longrightarrow Y$ be a map between
compacta. Then $f$ admits approximation by $k$-dimensional maps
if and only if there exists a map $g\colon X\longrightarrow {\bf
I}^k$ such that $\dim(f\bigtriangleup g)\le 0$.
\end{theorem}

There are a lot of papers devoted to Pasynkov's question
(\cite{pas},\cite{pasumn},\cite{passt},\cite{torun},\cite{stern},
\cite{lev},\cite{levlewis},\cite{tun_valov1}).
In \cite{pas} Pasynkov announced the following theorem to which the proof
appeared much later in \cite{pasumn} and  \cite{passt}.

\begin{theorem}
\label{Pasynkov} Let $f\colon X \longrightarrow Y$ be a
$k$-dimensional map between finite dimensional compacta. Then
there exists a map $g\colon X \longrightarrow {\bf I}^k$ such that
$\dim(f \bigtriangleup g)\leq 0.$
\end{theorem}

In \cite{torun} Torunczyk proved the following theorem which is
closely related to the theorem proved by Pasynkov.

\begin{theorem}
\label{Torunczyk} Let $f\colon X \longrightarrow Y$ be a
$k$-dimensional map between finite dimensional compacta. Then
there exists a $\sigma$-compact $A\subset X$ such that $\dim A
\leq k-1$ and $\dim f \mid_{X \setminus A} \leq 0.$
\end{theorem}

One can prove that for any map $f\colon X\longrightarrow Y$
between compacta the statement of theorem \ref{Pasynkov} holds
(for $f$) if and only if the statement of theorem \ref{Torunczyk}
holds (for $f$) (the proof can be found in \cite{lev}).

We improve the argument used by Torunzyk in \cite{torun} to prove
the next theorem and the implication
"$\Longleftarrow$" of theorem \ref{PU}.

\begin{theorem}
\label{2condit} Let $f\colon X \longrightarrow Y$ be a
$k$-dimensional map between $C$-compacta. Then

$(i) \  $ there exists a $\sigma$-compact subset $A \subset X$
such that \\$\dim A \leq k-1$ and
$\dim f \mid_{X \setminus A} \leq 0.$

$(ii) \  $ there exists a map $g\colon X \longrightarrow {\bf
I}^k$ such that $\dim (f \bigtriangleup g) \leq 0.$
\end{theorem}

\begin{remark}
After the paper had been already submitted for publication, the
author learned that using different techniques from those of this
paper V. Valov and H. Murat Tuncali \cite{tun_valov1} obtained a
more general result in the class of all metrizable spaces.
\end{remark}

In the next corollary by $e{-}\dim (X)$ we mean the extensional dimension
of a compact space $X$ introduced by A.N. Dranishnikov in \cite{dran}.

\begin{corollary}
Let $f\colon X\longrightarrow Y$ be a $k$-dimensional map between
$C$-compacta. Then $e{-}\dim (X) \leq e{-}\dim (Y\times {\bf
I}^k).$
\end{corollary}

\begin{proof}
From \cite{drusp} it follows that the extensional dimension
cannot be lowered by $0$-dimensional maps so the corollary is an
immediate consequence of theorem \ref{2condit}.
\end{proof}

One can understand the statement of the previous corollary
as a generalization of the classical Hurewicz formula.

\section{Proofs}
Further in this section we assume that every space $X$ is given
with a fixed metric $\rho_X$ which generates the same topology on
it. By $\rho_X (A,B)$ we mean the distance between subsets $A$ and
$B$ in the space $X$, namely, $\rho_X (A,B)=\inf \{\rho_X(a,b)|\
a\in A, b\in B \}$. The closure of a subset $A$ will be denoted
by $[A]$.

\begin{lemma}
\label{lemma1} Let $f\colon X\longrightarrow Y$ be a map between
compacta. Suppose that for any closed disjoint subsets $B$ and $C$
of $X$ there exists a closed subset $T$ of $X$ such that $\dim
T\le k-1$ and for any $y\in Y$ the set $T$ separates $f^{-1}(y)$
between $B$ and $C$. Then there exists a $\sigma$-compact subset
$A\subset X$ with $\dim A\le k-1$ such that $\dim f|_{X\setminus
A} \le 0$.
\end{lemma}

\begin{proof}
Take a countable open base ${\cal B}=\{U_\gamma|\  \gamma\in
\Gamma\}$ on $X$ such that the union of any finite number of sets
from ${\cal B}$ is again a member of ${\cal B}$. Define the set
$\Lambda\subset \Gamma\times\Gamma$ by the requirement:
$(\gamma,\mu)\in\Lambda$ if and only if $[U_\gamma]\cap [U_\mu]
=\varnothing$. Note that $\Lambda$ is countable. By assumption
for every pair $(\gamma,\mu)\in\Lambda$ there exists a set
$T_{(\gamma,\mu)}$ of dimension at most $k-1$ separating every
fiber $f^{-1}(y)$ between $[U_\gamma]$ and $[U_\mu]$. Now define
$A = \cup \{T_{(\gamma,\mu)}|\  (\gamma,\mu)\in\Lambda\}$. By
definition $A$ is  $\sigma$-compact and, obviously, $\dim A\le
k-1$. It is also easy to see that $\dim f\mid_{X\setminus A}\le
0$. Indeed, by the additivity property of the base ${\cal B}$ for
every pair of disjoint closed subsets $G$ and $H$ of a given
fiber $f^{-1}(y)$ there exists a pair $(\gamma,\mu)\in\Lambda$
such that $G\subset U_\gamma$ and $H\subset U_\mu$. Then
$T_{(\gamma,\mu)}\subset A$ separates $f^{-1}(y)$ between $G$ and
$H$. So, $\dim (f^{-1}(y)\setminus A)\le 0$.
\end{proof}

Let ${\cal F} = \mathbb N^0\cup \bigcup \{ \mathbb N^k : k \geq 1
\}$ be the union of all finite sequences of positive integers
plus empty sequence $\mathbb N^0=\{ *\}$. For every $i \in {\cal
F}$ let us denote by $|i|$ the length of the sequence $i$ and by
$(i,p)$ the sequence obtained by adding to $i$ a positive integer
$p$.

\begin{lemma}
\label{lemma2}
Let $f\colon X\longrightarrow Y$ be a map between
compacta. Let $B$ and $C$ be closed disjoint subsets of $X$.
Suppose that for every $i\in {\cal F}$ there are sets $U(i)$,
$V(i)$ and $F(i)$ such that:

(a) $F(i)$ is closed in $Y$, the sets $U(i)$ and $V(i)$ are open
subsets of $X$ and $[U(i)] \cap [V(i)] = \varnothing$;

(b) $U(*) \supset B, V(*) \supset C$ and $F(*) = Y$;

(c) $U(i,p)\supset U(i)\cap f^{-1}(F(i,p))$ and $V(i,p) \supset
V(i) \cap f^{-1}(F(i,p))$ for every $p\in\mathbb N$;

(d) $F(i)\subseteq\cup\{ F(i,p): p\in\mathbb N\}$ and ${\rm diam}
F(i) <\frac{1}{|i|}$;

(e) the set $E(i)=f^{-1}(F(i))\setminus (U(i)\cup V(i))$ admits
an open cover of order $k$ and diameter $\frac{1}{1+|i|}$;

(f) in notations of (e) the family $\{ E(i,p): p\in\mathbb N\}$
is discrete in $X$.
\\
Then there exists a closed subset $T$ of $X$ such that $\dim T\le
k-1$ and for any $y\in Y$ the set $T$ separates $f^{-1}(y)$
between $B$ and $C$.
\end{lemma}

\begin{proof}
We define the set $T$ in the following way:
$$T_n = \cup \{ E(i): |i| = n \} \ \ and \ \  T = \cap \{ T_n : n \geq 0 \}.$$
From property $(e)$ it follows that $\dim T\leq k-1$. Let us show
that for every $y\in Y$ the set $T$ separates $f^{-1}(y)$ between
$B$ and $C$. For every $y\in Y$ there exists a sequence $\{ i_n:
n\in\mathbb N\}$ such that
$$\{ y\}=F(i_1)\cap F(i_1,i_2)\cap ... \ .$$
Then $f^{-1}(y)\setminus T \subseteq U(y)\cup V(y)$ is a desired
partition. Here we denote by $U(y)$ and $V(y)$ the following sets
$$U(y) = f^{-1}(y) \cap \bigcup \{ U(i_1,...,i_p): p \in \mathbb N \},$$
$$V(y) = f^{-1}(y) \cap \bigcup \{ V(i_1,...,i_p): p \in \mathbb N \}.$$
\end{proof}

\begin{lemma}
\label{lemma3} Let $f\colon X\longrightarrow Y$ be a
$k$-dimensional map between $C$-compacta and $\epsilon$ be any
positive number. Let $U$ and $V$ be open subsets of $X$ with
$[U]\cap [V]=\emptyset$, and $F$ be a closed subset of $Y$ such
that $f(U)\cap f(V)\supseteq F$. Then there exist families of sets
$\{U_p\}$, $\{V_p\}$, and $\{F_p\}$ for $p\in \mathbb N$ such that:

(1) $F_p$ is closed in $Y$, the sets $U_p$ and $V_p$ are open subsets
of $X$ and $[U_p]\cap [V_p]=\varnothing$;

(2) $U_p\supseteq U\cap f^{-1}(F_p)$, $V_p\supseteq V\cap f^{-1}(F_p)$;

(3) $F\subseteq \bigcup \{ F_p: p\in \mathbb N \}$ and 
${\rm diam} F_p<\epsilon$;

(4) the set $E_p=f^{-1}(F_p)\setminus (U_p\cup V_p)$ admits an open cover
of order $k$ and diameter $\epsilon$;

(5) in notations of (4) the family $\{ E_p: p\in\mathbb N \}$ is discrete
in $X$.
\end{lemma}

\begin{proof}
Let $\{ W_l: l\in\mathbb N\}$ be a countable sequence of open
disjoint sets such that each of them separates $X$ between $[U]$
and $[V]$. For every $y\in F$ let $P_l(y)\subset W_l$ be a closed
$(k-1)$-dimensional set separating $f^{-1}(y)$ between $[U]$ and
$[V]$. Let $Q_l(y)\subset W_l$ be an open neighborhood of
$P_l(y)$ admitting a finite open cover of size $\epsilon$
and order $k$. As the map $f$ is closed there exists a
neighborhood $G_l(y)$ of $y$ in $F$ such that $f^{-1}([G_l(y)])
\cap Q_l(y)$ separates $f^{-1}([G_l(y)])$ between $[U]$ and $[V]$
and ${\rm diam} G_l(y)< \epsilon$. For every
$l\in\mathbb N$ the family $\alpha_l=\{ G_l(y) :y\in F\}$ is an
open cover of $F$. As $F$ is a $C$-compactum there exists a
finite sequence of finite disjoint open families of sets $\{
\mu_l :l\leq N\}$ such that each family $\mu_l$ refines the cover
$\alpha_l$ and $\mu=\cup\{ \mu_l: l\leq N\}$ is an open cover of
$F$. Further, for every $G\in\mu$ there are open subsets $U(G)$
and $V(G)$ of $X$ with disjoint closures such that
$$U(G)\supset f^{-1}(G)\cap [U], \ \ V(G)\supset f^{-1}(G)\cap [V]$$
and if $G$ is a member of $\mu_l$ then
$$f^{-1}(G)\setminus (U(G)\cup V(G))\subset Q_l(y)$$
for some $y\in F$. Let $\{ F(G): G\in\mu\}$ be a closed shrinking
of the cover $\mu$ and let
$$E(G)=f^{-1}(F(G))\setminus (U(G)\cup V(G))$$
for $G\in\mu$. Then for $G\in\mu_l$ and $H\in\mu_m$ we have
$E(G)\cap E(H)\subset W_l\cap W_m=\varnothing$. So the family $\{
E(G) : G\in\mu\}$ is discrete in $X$. Let us enumerate the
members of $\mu$: $\  G_1, G_2, ...$. To get the desired sets we
set
$$F_p=F(G_p), \  U_p=U(G_p), \  V_p=V(G_p).$$
\end{proof}

\begin{lemma}
\label{lemma4}
Let $f\colon X\longrightarrow Y$ be a
$k$-dimensional map between $C$-compacta. Then for any closed
disjoint subsets $B$ and $C$ of $X$ and for any $i\in {\cal F}$
there exist sets $U(i)$, $V(i)$ and $F(i)$ satisfying (a)--(f) of
Lemma \ref{lemma2}.
\end{lemma}

\begin{proof}
We will construct the sets $U(i)$, $V(i)$ and $F(i)$ by induction on
$|i|$. First set $F(*)=Y$ and $U(*)=U$, $V(*)=V$ for some open 
subsets $U$ and $V$ of $X$ with $[U]\cap [V]=\varnothing$. Assume 
the sets $U(i)$, $V(i)$ and $F(i)$ are already constructed and 
satisfy the conditions (a)-(f) of Lemma \ref{lemma2}. 
Now to get the sets  
$U(i,p)$, $V(i,p)$ and $F(i,p)$ for all $p\in\mathbb N$ apply
Lemma \ref{lemma3} to the sets $U=U(i)$, $V=V(i)$, $F=F(i)$ and
to $\epsilon=\frac{1}{1+|i|}$.
\end{proof}

\begin{lemma}
\label{lemma5} Let $f\colon X\longrightarrow Y$ be a map between
compacta admitting approximations by $k$-dimensional maps. Then
for any closed disjoint subsets $B$ and $C$ of $X$ there exist
sets $U(i)$, $V(i)$ and $F(i)$ satisfying (a)--(f) of Lemma
\ref{lemma2}.
\end{lemma}

\begin{proof}
The sets $U(i)$, $V(i)$ and $F(i)$ will be constructed by induction 
on $|i|$. First set $F(*)=Y$ and $U(*)=U$, $V(*)=V$ for some open 
subsets $U$ and $V$ of $X$ with $[U]\cap [V]=\varnothing$. Assume 
the sets $U(i)$, $V(i)$ and $F(i)$ are already constructed and 
satisfy the conditions (a)-(f) of Lemma \ref{lemma2}. Take
$\epsilon=\min \{ \frac{\rho(U(i),V(i))}{4}, \frac{1}{1+|i|} \}$.
By assumption, there exists a commutative diagram of the following 
form:

$$
\begin{array}{rcc}
 X &\stackrel{\kappa_X}{\longrightarrow}& |{\cal K}| \\
 f \downarrow &&\downarrow\lefteqn{p}\\
 Y &\stackrel{\kappa_Y}{\longrightarrow}& |{\cal L}|,
\end{array}
$$
\\
where $\kappa_X$ and $\kappa_Y$ are maps with $\epsilon$-small fibers.
Let $G=\kappa_X([U(i)])$, $H=\kappa_X([V(i)])$
and $F=\kappa_Y(F(i))$. Note that $G\cap H=\varnothing$.
Let $U$ and $V$ be open subsets of $|{\cal K}|$ with 
$U\supseteq G$, $V\supseteq H$ and $[U]\cap[V]=\varnothing$.
Let $\lambda_1$ be a Lebesgue number of some open covering on $|{\cal K}|$
whose preimage under the map $\kappa_X$ is an $\epsilon$-small 
covering on $X$.
Let $\lambda_2$ be a number defined similarly for $|{\cal L}|$ and $\kappa_Y$.
Let $\lambda=\min\{\lambda_1,\lambda_2\}$. 
Apply Lemma \ref{lemma3} to the sets $U$, $V$, $F$ and to $\lambda$
to produce the sets $U_p$, $V_p$ and $F_p$ for all $p\in \mathbb N$
satisfying conditions (1)-(5) of Lemma \ref{lemma3}. 
Now set $U(i,p)=\kappa_X^{-1}(U_p)$, $V(i,p)=\kappa_X^{-1}(V_p)$ and
$F(i,p)=\kappa_Y^{-1}(F_p)$. Since taking a preimage preserves
intersections, unions and subtractions, the sets $U(i,p)$, $V(i,p)$
and $F(i,p)$ satisfy the conditions (a)-(f).
\end{proof}

\begin{proof}[Proof of theorem \ref{2condit}]
The statements $(i)$ and $(ii)$ are equivalent (\cite{lev}), so,
it is sufficient to prove only $(i)$. But $(i)$ immediately 
follows from Lemmas \ref{lemma1}, \ref{lemma2} and \ref{lemma4}.
\end{proof}

Let $\gamma$ be an open cover on $X$. Then by $N_\gamma$ we mean
the nerve of the cover $\gamma$. Let $\{ a_\alpha : \alpha\in
A\}$ be some partition of unity on $X$ subordinated to the
locally finite cover $\gamma$. Then the canonical map defined by
the partition of unity $\{ a_\alpha : \alpha\in A \}$ is a map
$\kappa\colon X\longrightarrow |N_\gamma|$ defined by the
following formula
$$\kappa(x) = \sum_{\alpha \in A} a_\alpha(x)\cdot \alpha.$$
If $\tau$ is some triangulation of the polyhedron $P$, then
by $St(a,\tau)$ we mean the star of the vertex
$a\in\tau$ with respect to triangulation $\tau$, {\sl i.e.} the union of all
open simplices having $a$ as a vertex.

\begin{proof}[Proof of theorem \ref{PU}]
Let $f\colon X\longrightarrow Y$ be a map between compacta
admitting approximation by $k$-dimensional simplicial maps. 
By \cite{lev} to show that there exists a map 
$g\colon X\longrightarrow {\bf I}^k$ with $\dim(f\bigtriangleup g)\le 0$
it is safficient to find a $\sigma$-compact subset $A$ in $X$ 
of dimension at most $k-1$ such that $\dim f\mid_{X\setminus
A}\le 0$. The existence of such subset $A$ follows immediately from 
Lemmas \ref{lemma1}, \ref{lemma2} and \ref{lemma5}.

Now suppose there exists a map $g\colon X\longrightarrow {\bf
I}^k$ such that $\dim (f \bigtriangleup g) \leq 0.$ For every
$(y,t)\in Y\times {\bf I}^k$ there exists a finite disjoint
family of open sets $\nu_{(y,t)} = \{ V_\gamma :
\gamma\in\Gamma_{(y,t)}\}$ such that $(f\bigtriangleup
g)^{-1}(y,t)\subset\cup\nu_{(y,t)}$ and
$\nu_{(y,t)}\succ\omega_X.$ Let $O_{(y,t)}$ be an open
neighborhood of $(y,t)$ in $Y\times {\bf I}^k$ such that
$(f\bigtriangleup g)^{-1}O_{(y,t)}\subset\cup\nu_{(y,t)}.$ Let
$\varsigma=\{ U_\alpha : \alpha\in A \}$ and $\iota=\{ I_\delta :
\delta\in D\}$ be finite open covers of the spaces $Y$ and ${\bf
I}^k$ such that:

(a$_1$) $\varsigma\succ\omega_Y$;

(b$_1$) the order of $\varsigma$ does not exceed $\dim Y+1$;

(c$_1$) $(\varsigma\times\iota) =
\{ U_\alpha\times I_\delta : (\alpha,\delta)\in A\times D\}\succ
\{ O_{(y,t)} : (y,t)\in Y\times {\bf I}^k\}$.

The partition of unity $\{ u_\alpha : \alpha\in A\}$ on $Y$
subordinated to the cover $\varsigma$ gives rise to the canonical
map $\mu\colon Y\longrightarrow |{\cal N}_\varsigma|$. Then the
map $\mu\times id\colon Y\times {\bf I}^k\longrightarrow |{\cal
N}_\varsigma|\times {\bf I}^k$ is an $(\varsigma\times\iota)$-map.
By $\pi\colon |{\cal N}_\varsigma|\times {\bf I}^k\longrightarrow
|{\cal N}_\varsigma|$ we denote the projection. Let $\tau$ and
$\theta$ be such triangulations on polyhedra $|{\cal
N}_\varsigma|\times {\bf I}^k$ and $|{\cal N}_\varsigma|$
respectively such that the following conditions are satisfied:

(1$^\prime$) $\pi\colon |{\cal N}_\varsigma|\times {\bf
I}^k\longrightarrow |{\cal N}_\varsigma|$ is a simplicial map
relative to the triangulations $\tau$ and $\theta$;

(2$^\prime$)
$\{ (\mu\times id)^{-1}(St(a,\tau)) :a\in \tau \} \succ \varsigma\times\iota$;

(3$^\prime$) $\{ \mu^{-1}(St(b,\theta)) :b\in \theta \} \succ \varsigma$.

Let us define $\xi=\{ St(a,\tau) :a\in\tau\}$ and $\zeta=\{
St(b,\theta) :b\in \theta \}$. Define the partition of unity $\{
w_a: a\in\tau\}$ on $|{\cal N}_\varsigma|\times {\bf I}^k$
subordinated to the cover $\xi$ by letting $w_a(z)$ to be the
barycentric coordinate of $z\in |{\cal N}_\varsigma|\times {\bf
I}^k$ with respect to the vertex $a\in\tau$. Analogously, define
the partition of unity $\{ v_b: b\in\theta\}$ on $|{\cal
N}_\varsigma|$ subordinated to the cover $\zeta$. Note that the
projection $\pi\colon |{\cal N}_\varsigma|\times {\bf
I}^k\longrightarrow |{\cal N}_\varsigma|$ sends the stars of the
vertices of the triangulation $\tau$ to the stars of the vertices
of the triangulation $\theta$. That is why there is a simplicial
map $\varpi\colon {\cal N}_\xi\longrightarrow {\cal N}_\zeta$
between the nerves of the covers $\xi$ and $\zeta$. Moreover, the
following diagram commutes.

$$
\begin{array}{ccc}
 |{\cal N}_\varsigma|\times {\bf I}^k &\stackrel{\psi}{\longrightarrow}& |{\cal N}_\xi| \\
 \pi\downarrow &&\downarrow\lefteqn{\varpi}\\
 |{\cal N}_\varsigma| &\stackrel{\phi}{\longrightarrow}& |{\cal N}_\zeta|.
\end{array}
$$
\\
Here $\psi$ and $\phi$ are canonical maps defined by the partitions
of unity $\{ w_a: a\in\tau \}$ and $\{ v_b: b\in \theta \}$ respectively.
Let us remark that $\dim \varpi\le k$.
By $W_a$ we will denote the set $(\mu\times id)^{-1}(St(a,\tau))$,
by $\lambda$ the cover $\{ W_a: a\in\tau \}$ and $\eta$ is $\mu^{-1}(\zeta)$.
Further, we set $w^*_a = w_a\circ (\mu\times id)$ for each $a\in\tau$ and
$v^*_b = v_b\circ \mu$ for each $b\in\theta$.
The partitions of unity $\{ w^*_a : a\in \tau \}$ on $Y\times {\bf I}^k$ and
$\{ v^*_b : b\in \theta \}$ on $Y$ are subordinated to the covers
$\lambda$ and $\eta$ respectively.
We set ${\cal K}^\prime={\cal N}_\lambda$ and
${\cal L}={\cal N}_\eta$. The simplicial complexes
${\cal K}^\prime$ and ${\cal L}$ are isomorphic to ${\cal N}_\xi$ and
${\cal N}_\zeta$ that is why the simplicial map
$q:|{\cal K}^\prime|\longrightarrow |{\cal L}|$ is defined and
the following diagram commutes.

$$
\begin{array}{ccc}
 \quad\enskip Y\times {\bf I}^k &\stackrel{\varphi}{\longrightarrow}& |{\cal K}^\prime| \\
 pr\downarrow &&\downarrow\lefteqn{q}\\
 \quad\enskip Y &\stackrel{\kappa_Y}{\longrightarrow}& |{\cal L}|.
\end{array}
$$
\\
Here $\varphi$ and $\kappa_Y$ are canonical maps defined by the
partitions of unity
$\{ w^*_a : a\in\tau\}$ and $\{ v^*_b : b\in\theta\}$ respectively.
Obviously, $\dim q\le k$.

Recall that $\varphi$ is an $\{ O_{(y,t)}\}$-map.
For every $a\in\tau$ pick a point $(y,t)_a$ such that
$W_a\subset O_{(y,t)_a}$.
Let $w^{**}_a=w^*_a\circ (f\bigtriangleup g)$
and $B_a=\Gamma_{(y,t)_a}$. Then
$$supp \  (w^{**}_a) \subset
(f \bigtriangleup g)^{-1}(W_a)\subset (f \bigtriangleup g)^{-1}(O_{(y,t)_a}).$$

Consequently, $\cup\nu_{(y,t)_a}\supset supp (w^{**}_a)$. As the
family $\nu_{(y,t)_a}$ is disjoint, there exists a family of
non-negative functions $\{ b_\beta : \beta\in B_a\}$ such that
$w^{**}_a =\sum_{\beta\in B_a} b_\beta$ and $supp
(b_\beta)\subset V_\beta$. Let $B=\cup\{ B_a : a\in\tau\}$ and
${\cal V}=\{ V_\beta\cap (f \bigtriangleup g)^{-1}(W_a) :
a\in\tau,\beta\in B_a\}$. The family $\{ b_\beta : \beta\in B\}$
is a partition of unity on $X$ subordinated to the cover ${\cal
V}$. Let ${\cal K}$ be the nerve of the cover ${\cal V}$ and
$\kappa_X\colon X \longrightarrow |{\cal K}|$ the canonical map
defined by the partition of unity $\{ b_\beta: \beta\in B \}$. We
define a simplicial map $p'\colon |{\cal K}|\longrightarrow |{\cal
K'}|$ by requiring that the vertex $\beta \in B_a$ goes to $a$.
Clearly, the map $p^\prime$ is finite-to-one. Indeed, no two
vertices in $B_a$ are connected by an edge. That is why the
restriction of $p^\prime$ to any simplex is a homeomorphism.
Finally we define the desired $k$-dimensional simplicial map
$p\colon |{\cal K}|\longrightarrow |{\cal L}|$ as the composition
$p = q \circ p'$. Moreover by (b$_1$) we have $\dim |{\cal L}|\le
\dim Y$ and $\dim |{\cal K}|\le \dim Y+k$ since $p$ is
$k$-dimensional.
\end{proof}

\begin{proof}[Proof of theorem \ref{th_C-comp}]
The theorem is an immediate consequence of theorems \ref{PU} and
\ref{2condit} and the remark at the end of the proof of theorem \ref{PU}.
\end{proof}

The following results of M.Levin and Y.Sternfeld are needed to prove theorems
\ref{arbcomp} and \ref{bingcomp}:

\begin{theorem}[\cite{lev}]
\label{th_lev} Let $f\colon X\longrightarrow Y$ be a
$k$-dimensional map between compacta. Then there exists a map
$g\colon X\longrightarrow {\bf I}^{k+1}$ such that $\dim(f
\bigtriangleup g) \leq 0.$
\end{theorem}

\begin{theorem}[\cite{stern}]
\label{th_stern} Let $f\colon X\longrightarrow Y$ be a
$k$-dimensional map of Bing compacta. Then there exists a map
$g\colon X\longrightarrow {\bf I}^k$ such that $\dim(f
\bigtriangleup g) \leq 0.$
\end{theorem}

\begin{proof}[Proof of theorem \ref{arbcomp}]
The theorem is an immediate consequence of theorems
\ref{th_lev} and \ref{PU}.
\end{proof}

\begin{proof}[Proof of theorem \ref{bingcomp}]
The theorem is an immediate consequence of theorems
\ref{th_stern} and \ref{PU}.
\end{proof}

\section*{Acknowledgements}

I would like to thank A. N. Dranishnikov, S. A. Melikhov, B.A. Pasynkov,
R. R. Sadykov and V. V. Uspenskij for useful discussions and important 
remarks. I also would like to thank the unknown referee for many important 
remarks concerning the structure of the paper.

\section*{}
{\small
Moscow State University and University of Florida\\
E-mail address: turygin@math.ufl.edu}

\begin{thebibliography}{99}
\bibitem{dran} A. N. Dranishnikov,
{\sl The Eilenberg-Borsuk theorem for maps into arbitrary complexes,}
Math. Sbornik vol. 185(1994), no. 4, 81-90.
\bibitem{drusp} A. N. Dranishnikov {\sl and} V. V. Uspenskij,
{\sl Light maps and extensional dimension,} Topology Appl., 80 (1997) 91-99.
\bibitem{lev} M. Levin, {\sl Bing maps and finite dimensional maps,}
Fund. Math. 151(1) (1996) 47-52.
\bibitem{levlewis} M. Levin {\sl and} W. Lewis,
{\sl Some mapping theorems for extensional dimension}, arXiv:math.GN/0103199
\bibitem{pas} B. A. Pasynkov, {\sl The dimension and geometry of mappings,}
(Russian) Doklad. Akad. Nauk. SSSR 221(1975), 543-546.
\bibitem{pasumn} B. A. Pasynkov, {\sl Theorem on $\omega$-mappings for
mappings,} (Russian) Uspehi. Math. Nauk. vol.39 (1984), no.
5(239), 107-130.
\bibitem{passt} B. A. Pasynkov,
{\sl On the geometry of continuous mappings of finite dimensional
metrizable compacta,} Proc. Steklov Inst. Math. vol. 212(1996),
138-162.
\bibitem{stern} Y. Sternfeld, {\sl On finite-dimensional maps and other maps
with "small" fibers,} Fund. Math., 147 (1995), 127-133.
\bibitem{torun} H. Torunczyk, {\sl Finite to one restrictions of continuous
functions,} Fund. Math. 125 (1985), 237-249.
\bibitem{tun_valov1} H. Murat Tuncali, V. Valov,
{\sl On dimensionally restricted maps,} Fund. Math. 175
(1)(2002), 35--52.
\bibitem{usp} V. V. Uspenskij, {\sl A selection theorem for $C$-spases,}
Topology Appl. 85 (1998) 351-374.
\end{thebibliography}
\end{document}